# On continuous-time autoregressive fractionally integrated moving average processes

HENGHSIU TSAI

*Institute of Statistical Science, Academia Sinica, Taipei, Taiwan 115, Republic of China.
E-mail: htsai@stat.sinica.edu.tw*

In this paper, we consider a continuous-time autoregressive fractionally integrated moving average (CARFIMA) model, which is defined as the stationary solution of a stochastic differential equation driven by a standard fractional Brownian motion. Like the discrete-time ARFIMA model, the CARFIMA model is useful for studying time series with short memory, long memory and antipersistence. We investigate the stationarity of the model and derive its covariance structure. In addition, we derive the spectral density function of a stationary CARFIMA process.

## 1. Introduction

It is well known that the long-range dependence properties of time series data have diverse applications in many fields, including hydrology, economics and telecommunications (see Bloomfield [8], Sowell [35], Beran [4], Robinson [32], Baillie [3] and Ray and Tsay [31]). Autoregressive fractionally integrated moving average (ARFIMA) models are a popular class of discrete-time long memory processes (see Granger and Joyeux [19] and Hosking [21]). For continuous-time long memory modeling, see Viano *et al.* [40], Chambers [13], Comte [14], Comte and Renault [15], Brockwell and Marquardt [12], Tsai and Chan [37] and Tsai and Chan [38].

In contrast to long-range dependence, which corresponds to the singularity of the spectrum at the origin, antipersistent time series are covariance-stationary processes with zero spectral density at the origin (Dittmann and Granger [16], Bondon and Palma [9]). Interesting examples of antipersistent processes include: Kolmogorov [24], which models the local structure of turbulence in incompressible viscous fluids; Ausloos and Ivanova [2], which models the temporal correlations of fluctuations in the Southern Oscillation index (SOI) signal; Simonsen [34], which measures correlations in the Nordic electricity spot market. For financial applications of antipersistent processes, see Peters [29] and Shiryaev [33], both of which model the implied and realized volatility of the S & P 500 index. For further examples of antipersistence modeling, see Otway [28], Beran and Mazzola [7] and Metzler [26].

Beran and Feng [5] observed that antipersistent processes can be generated by overdifferencing. In other words, a slightly non-stationary process lying between a stationary long-memory process and a random walk becomes antipersistent after the first differencing. For example, Beran *et al.* [6] found that the first-differenced daily world copper price from January 2, 1997 to September 2, 1998 comprised an antipersistent process and a non-zero mean function. This implies that the original process was a trend function, plus a nonstationary error process whose first difference was antipersistent. Karuppiah and Los [23] also noted that many intra-day foreign exchange rate series are antipersistent.

Recently, Tsai and Chan [37] proposed the continuous-time autoregressive fractionally integrated moving average (CARFIMA) model, which is defined as the stationary solution of a stochastic differential equation driven by standard fractional Brownian motion with the Hurst parameter $1/2 < H < 1$. The model is useful for studying time series data that exhibits long-range dependence properties. In [37], Tsai and Chan derived the autocovariance function of the stationary CARFIMA process and in [38], they derived its spectral density function; both functions are derived under the condition $1/2 < H < 1$.

In this paper, we focus on the stationarity, autocovariance function, spectral density function and other properties of the CARFIMA process with the Hurst parameter $0 < H < 1$. The CARFIMA process with $0 < H < 1/2$ is an antipersistent process that is useful for studying time series data exhibiting antipersistent properties. For $H = 1/2$, the CARFIMA process becomes a continuous-time autoregressive moving average (CARMA) process. We present the main results in Section 2, followed by their proofs in Section 3. Finally, in Section 4, we make some concluding remarks.

## 2. Main results

Heuristically, a CARFIMA($p$, H, $q$) process $\{Y_t\}$ is defined as the solution of a $p$th order stochastic differential equation with suitable initial conditions. It is driven by a standard fractional Brownian motion with the Hurst parameter $H$ and its derivatives up to and including the order $0 \leq q < p$. Specifically, for $t \geq 0$,

$$Y_t^{(p)} - \alpha_p Y_t^{(p-1)} - \cdots - \alpha_1 Y_t - \alpha_0 \qquad (1)$$
$$= \sigma\{B_{t,H}^{(1)} + \beta_1 B_{t,H}^{(2)} + \cdots + \beta_q B_{t,H}^{(q+1)}\},$$

where $\{B_{t,H} = B_t^H, t \geq 0\}$ is a standard fractional Brownian motion with the Hurst parameter $0 < H < 1$, $dY_t^{(j-1)} = Y_t^{(j)} dt, j = 1, \ldots, p-1$. The superscript $^{(j)}$ denotes a $j$-fold differentiation with respect to $t$. We assume that $\sigma > 0$, $\alpha_1 \neq 0$ and $\beta_q \neq 0$.

Let $0 < H < 1$ be a fixed number. It is well known (see, e.g., Duncan *et al.* [17]) that there exists a Gaussian stochastic process $\{B_t^H, t \geq 0\}$, which satisfies the following three properties:

(a) its initial condition is $B_0^H = 0$;

(b) it is of zero mean, that is, $\mathrm{E}(B_t^H) = 0$ for all $t \geq 0$;
(c) its covariance kernel is defined as

$$\mathrm{E}(B_t^H B_s^H) = \tfrac{1}{2}\{|t|^{2H} + |s|^{2H} - |t-s|^{2H}\}$$

for all $s, t \geq 0$. The Gaussian process $\{B_t^H\}$ is called a standard fractional Brownian motion with the Hurst parameter $H$. A standard fractional Brownian motion with $H = 1/2$ is the same as a standard Brownian motion. In addition, a fractional Brownian motion defined over non-negative $t$ can be extended and defined for all real numbers. Specifically, for $H \in (0,1)$, we have, for $-\infty < t < \infty$,

$$B_t^H = D_H \int_{-\infty}^{\infty} \{(t-u)_+^{H-1/2} - (-u)_+^{H-1/2}\}\,\mathrm{d}B_u,$$

where $D_H = [2H\Gamma(3/2-H)/\{\Gamma(H+1/2)\Gamma(2-2H)\}]^{1/2}$, $\Gamma(\cdot)$ is the Gamma function and $\{B_u\}$ is a *two-sided* standard Brownian motion with covariance $\mathrm{E}(B_s B_t) = \min(|s|,|t|)$ if $s$ and $t$ have the same sign and $0$ if $s$ and $t$ have different signs (see Taqqu [36]).

Let $Z$ be a set of integers. The stationary process $\{F_t\}_{t \in Z}$, defined by $F_t = B_{t+1}^H - B_t^H$, is called *fractional Gaussian noise*. The autocovariance of the noise is

$$\gamma_F(k) = \tfrac{1}{2}\{|k+1|^{2H} - 2|k|^{2H} + |k-1|^{2H}\}.$$

It can be easily demonstrated (see, e.g., Taqqu [36]) that for $k \neq 0$, $\gamma_F(k) = 0$ if $H = 1/2$, $\gamma_F(k) < 0$ if $0 < H < 1/2$ and $\gamma_F(k) > 0$ if $1/2 < H < 1$. For $H \neq 1/2$,

$$\gamma_F(k) \sim H(2H-1)|k|^{2H-2} \qquad \text{as } k \to \infty.$$

For $1/2 < H < 1$, $\sum_{k=-\infty}^{\infty} \gamma_F(k) = \infty$, so the noise is said to be persistent, to have long memory or to be long-range dependent. For $0 < H < 1/2$, $\sum_{k=-\infty}^{\infty} \gamma_F(k) = 0$; therefore, the process is negatively autocorrelated at all positive lags and the noise is said to be antipersistent. For $H = 1/2$, all correlations of the process at non-zero lags are zero, that is, $F_t$ and $F_s$ are uncorrelated for $t \neq s$. The spectral density of the noise is given by (see Beran [4])

$$f_F(\omega) = \{2(1-\cos\omega)\} \sum_{k=-\infty}^{\infty} |\omega + 2k\pi|^{-2H-1}.$$

Note that $f_F(\omega)$ is $O(|\omega|^{1-2H})$ for $\omega \to 0$.

Taqqu [36] considers continuous-time fractional Gaussian noise $\{C_t\}_{t \in R}$, where $R$ is a set of real numbers, namely $C_t = B_{t+1}^H - B_t^H$, and shows that the spectral density of continuous-time fractional Gaussian noise is given by

$$f_C(\omega) = \frac{1}{2\pi}\sigma^2 \Gamma(2H+1) \sin(\pi H) |\mathrm{e}^{\mathrm{i}\omega} - 1|^2 |\omega|^{-2H-1}$$

$$= \frac{2}{\pi}\sigma^2 \Gamma(2H+1) \sin(\pi H) \sin^2\left(\frac{\omega}{2}\right) |\omega|^{-2H-1}, \qquad \omega \in R.$$

Fractional Brownian motion is not differentiable (Mandelbrot and Van Ness [25]), so the stochastic equation (1) has to be appropriately interpreted as an integral equation, as explained below. Analogous to continuous-time ARMA processes (see, e.g., Brockwell [10]), equation (1) can be equivalently cast in terms of the *observation* and *state* equations

$$Y_t = \beta' X_t, \qquad t \geq 0, \tag{2}$$

$$dX_t = (AX_t + \alpha_0 \delta_p)\,dt + \sigma \delta_p \,dB_t^H, \tag{3}$$

where the prime superscript denotes taking the transpose,

$$A = \begin{bmatrix} 0 & 1 & 0 & \cdots & 0 \\ 0 & 0 & 1 & \cdots & 0 \\ \vdots & \vdots & \vdots & \ddots & \vdots \\ 0 & 0 & 0 & \cdots & 1 \\ \alpha_1 & \alpha_2 & \alpha_3 & \cdots & \alpha_p \end{bmatrix}, \qquad X_t = \begin{bmatrix} X_t^{(0)} \\ X_t^{(1)} \\ \vdots \\ X_t^{(p-2)} \\ X_t^{(p-1)} \end{bmatrix},$$

$$\delta_p = \begin{bmatrix} 0 \\ 0 \\ \vdots \\ 0 \\ 1 \end{bmatrix}, \qquad \beta = \begin{bmatrix} 1 \\ \beta_1 \\ \vdots \\ \beta_{p-2} \\ \beta_{p-1} \end{bmatrix}$$

and $\beta_j = 0$ for $j > q$. See Tsai and Chan [39] for a proof of the equivalence between equation (1) and the set of equations (2) and (3). The process $\{Y_t, t \geq 0\}$ is said to be a CARFIMA($p$, H, $q$) process with the parameter $(\theta, \sigma) = (\alpha_0, \ldots, \alpha_p, \beta_1, \ldots, \beta_q, H, \sigma)$ if $Y_t = \beta' X_t$, where $X_t$ is the solution of (3) for the initial condition $X_0$. Similarly to equation (5) of Tsai and Chan [37], the solution $\{X_t, t \geq 0\}$ of (3) can be written as

$$X_t = e^{At} X_0 + \alpha_0 \int_0^t e^{A(t-u)} \delta_p \,du + \sigma \int_0^t e^{A(t-u)} \delta_p \,dB_u^H, \tag{4}$$

where

$$e^{At} = I_p + \sum_{n=1}^{\infty} \{(At)^n (n!)^{-1}\}$$

and $I_p$ is the $p \times p$ identity matrix. The stochastic integration in (4) is defined in terms of the limit of Riemann sums because it only involves deterministic integrands.

For a random initial condition $X_0$, the mean vector of $\{X_t\}$, denoted by $\mu_{X,t}$, satisfies the following equation:

$$\mu_{X,t} = e^{At} \mu_{X,0} + \frac{\alpha_0}{\alpha_1}(e^{At} - I_p)\delta_1, \tag{5}$$

where $\delta_1 = [1, 0, \ldots, 0]'$. If $\mu_{X,0}$ is chosen to be $-(\alpha_0/\alpha_1)\delta_1$, then $\mu_{X,t}$ becomes $-(\alpha_0/\alpha_1)\delta_1$, which is independent of $t$. If all the eigenvalues of $A$ have negative real parts, then it can be easily shown that (2) and (4) imply that, when $t \to \infty$, $Y_t$ converges in distribution to a normal random variable with mean $-\alpha_0/\alpha_1$ and variance $V_Y$, where $V_Y = \gamma_Y(0)$ is defined in equation (11). Thus, the stationary solution, if it exists, must be Gaussian.

The stationary CARFIMA process defined over non-negative $t$ can be extended so that it is a stationary process over all real $t$. For simplicity, we assume that $\alpha_0 = 0$. Then, provided the eigenvalues of $A$ all have negative real parts, we can show that the process $\{X_t\}$ defined by

$$X_t = \sigma \int_{-\infty}^{t} e^{A(t-u)} \delta_p \, dB_u^H \tag{6}$$

is a strictly stationary solution of (3) for $t \in (-\infty, \infty)$ with the corresponding CARFIMA process

$$Y_t = \sigma \int_{-\infty}^{t} \beta' e^{A(t-u)} \delta_p \, dB_u^H. \tag{7}$$

The proof of (6) is similar to that of (4) and is hence omitted.

For $1/2 < H < 1$ and $f, g \in L^2(R; R) \cap L^1(R; R)$, Gripenberg and Norros [20] proved that

$$\operatorname{cov}\left(\int_{-\infty}^{\infty} f(u) \, dB_u^H, \int_{-\infty}^{\infty} g(v) \, dB_v^H\right)$$
$$= H(2H-1) \int_{-\infty}^{\infty} \int_{-\infty}^{\infty} f(u)g(v)|u-v|^{2H-2} \, du \, dv. \tag{8}$$

For $H < 1/2$, Norros et al. [27] defines a class of stochastic integrals with the fractional Brownian integrator and deterministic integrands that are functions of bounded variation. However, they only consider integrals over the interval $[0, \infty)$. In our case, integrals over an interval of the form $(-\infty, t]$ for finite $t$ are needed. For the definition of such stochastic integrals, we outline an approach similar to Norros et al. [27]. Let $t$ be a fixed but arbitrary finite real number. Let $\Gamma$ denote the integral operator mapping a bounded-variation function $f(s), s \leq t$, to a function $\Gamma f(s), s \leq t$, defined by the equation

$$\Gamma f(s) = H f(t) |t-s|^{2H-1} \operatorname{sgn}(t-s)$$
$$+ H \int_{-\infty}^{t} |u-s|^{2H-1} \operatorname{sgn}(s-u) \, df(u),$$

where $\operatorname{sgn}(u)$ equals $-1$, 0 or 1, depending on whether $u$ is negative, zero or positive. Next, define an inner product between two functions $f$ and $g$ by the formula

$$\langle f, g \rangle_\Gamma = \int_{-\infty}^{t} g(v) \Gamma f(v) \, dv.$$

Define $L^2_\Gamma(-\infty, t]$ as the space of equivalence classes of measurable, bounded-variation functions $f$ such that $\langle f, f \rangle_\Gamma < \infty$. Consider the association of the simple function $1_{[a,b]}$ to $B^H_b - B^H_a$ that preserves the inner product, which follows from equation (9) below, and hence the association can be extended to an isometry between the Gaussian space spanned by $B^H_u, u \leq t$, and the function space $L^2_\Gamma(-\infty, t]$ so that the integral $\int_{-\infty}^t f(u) \, dB^H_u$ can be defined as the image of $f$ in the isometry. Consequently, provided that the integrals of the right-hand side of equation (9) exist, then for $0 < H < 1/2$ and all $s, t \in R$,

$$\text{cov}\left(\int_{-\infty}^s f(u) \, dB^H_u, \int_{-\infty}^t g(v) \, dB^H_v\right)$$

$$= H f(s) \int_{-\infty}^t |s-v|^{2H-1} \text{sgn}(s-v) g(v) \, dv$$

$$+ H \int_{-\infty}^t \int_{-\infty}^s g(v) |u-v|^{2H-1} \text{sgn}(v-u) \, df(u) \, dv. \tag{9}$$

To demonstrate the above formula, it suffices to consider the case $t = s$. The proof of equation (9) for simple functions is given in Section 3 and its validity for general functions then follows from the isometry alluded to above.

From the above discussion, we can derive the stationarity condition and the autocovariance function of the CARFIMA process with $0 < H < 1$. We state the stationarity condition of the CARFIMA process in Theorem 1.

**Theorem 1.** *Let $0 < H < 1$. Equation (4) with a deterministic initial condition admits an asymptotically stationary solution if and only if all the eigenvalues of $A$ have negative real parts. Moreover, under the preceding eigenvalue condition of $A$ and assuming the solution is stationary, $Y_0$ and $\{B^H_t, t \geq 0\}$ are jointly Gaussian with the covariances given by*

$$\text{cov}(Y_0, B^H_t) = H \sigma \beta' \int_0^\infty e^{Au} \delta_p \{(u+t)^{2H-1} - u^{2H-1}\} \, du, \tag{10}$$

*and the stationary mean of $\{Y_t\}$ equalling $\mu_Y = -\alpha_0/\alpha_1$.*

For a random initial condition $Y_0$, which may be correlated with the fractional Brownian innovation process, it can be verified that the sufficiency part of Theorem 1 continues to hold if $Y_0$ has finite variance. Furthermore, the theorem implies that, under stationarity, $Y_0$ and the fractional Brownian innovations $B^H_t, t \geq 0$, are generally correlated when $H \neq 1/2$. This contrasts with the case $H = 1/2$, where the stationary distribution of $Y_0$ is independent of the standard Brownian motion.

In Theorem 2, we use the covariance formulae of stochastic integrals given in expressions (8) and (9) to calculate the autocovariance function of $\{Y_t\}$. In part (a), the autocovariance function of the CARFIMA process is expressed in terms of three integrals. In part (b), the eigenvalues of the companion matrix $A$ are distinct, so we have a closed form of the autocovariance function. Then, (c) describes the asymptotic expression for the autocovariance function with $H \neq 1/2$.

**Theorem 2.** *Let $0 < H < 1$.*
*(a) Under stationarity, for $h \geq 0$, the autocovariance function of $\{Y_t\}$ equals*

$$\begin{aligned}
\gamma_Y(h) &:= \operatorname{cov}(Y_{t+h}, Y_t) \\
&= H\beta' A \mathrm{e}^{Ah}\left(\int_0^h \mathrm{e}^{-Au} u^{2H-1}\,\mathrm{d}u\right) V^*\beta \\
&\quad - H\beta' A \mathrm{e}^{-Ah}\left(\int_h^\infty \mathrm{e}^{Au} u^{2H-1}\,\mathrm{d}u\right) V^*\beta - H\beta' A \mathrm{e}^{Ah}\left(\int_0^\infty \mathrm{e}^{Au} u^{2H-1}\,\mathrm{d}u\right) V^*\beta,
\end{aligned} \tag{11}$$

*where $V^* = \sigma^2 \int_0^\infty \mathrm{e}^{Au} \delta_p \delta'_p \mathrm{e}^{A'u}\,\mathrm{d}u$.*
*(b) Under stationarity, and when the eigenvalues $\lambda_1, \ldots, \lambda_p$ of the companion matrix $A$ are distinct,*

$$\gamma_Y(h) = \frac{\sigma^2}{2}\Gamma(2H+1)\sum_{i=1}^p \frac{\beta(\lambda_i)\beta(-\lambda_i)}{\alpha^{(1)}(\lambda_i)\alpha(-\lambda_i)} u(H, \lambda_i, h), \tag{12}$$

*where $h \geq 0$, $\alpha(z) = z^p - \alpha_p z^{p-1} - \cdots - \alpha_1$, $\alpha^{(1)}(\cdot)$ denotes its first derivative, $\beta(z) = 1 + \beta_1 z + \beta_2 z^2 + \cdots + \beta_q z^q$,*

$$\begin{aligned}
u(H, \lambda, h) &= 2(-\lambda)^{1-2H}\cosh(\lambda h) + \lambda^{1-2H}\mathrm{e}^{\lambda h}P(2H, \lambda h) \\
&\quad - (-\lambda)^{1-2H}\mathrm{e}^{-\lambda h}P(2H, -\lambda h)
\end{aligned} \tag{13}$$

*and $P(a, z) = \int_0^z \mathrm{e}^{-u} u^{a-1}\,\mathrm{d}u/\Gamma(a)$, where the integration is along the radial line in the complex plane from 0 to $z$.*
*(c) For $H \neq 1/2$, as $h \to \infty$, we have the asymptotic expression*

$$\gamma_Y(h) \sim \sigma^2 H(2H-1)\frac{\beta^2(0)}{\alpha^2(0)} h^{2H-2}, \tag{14}$$

*where "$\sim$" means that the ratio of the left- and right-hand sides converges to 1.*

For $H = 1/2$, it is easy to verify that Theorem 2(a) and (b) can be simplified to

$$\begin{aligned}
\gamma_Y(h) &= \beta' \mathrm{e}^{Ah} V^* \beta \\
&= \sigma^2 \sum_{i=1}^p \frac{\beta(\lambda_i)\beta(-\lambda_i)}{\alpha^{(1)}(\lambda_i)\alpha(-\lambda_i)} \mathrm{e}^{\lambda_i h},
\end{aligned}$$

which is consistent with the autocovariance function of the short-memory CARMA model proposed by Brockwell [11]. Let $\Gamma(a, z) = \int_z^\infty \mathrm{e}^{-u} u^{a-1}\,\mathrm{d}u$ be an incomplete Gamma function with complex arguments. Then, the fact that $\Gamma(a+1, z) = a\Gamma(a, z) + z^a \mathrm{e}^{-z}$ can be used to show that, for $1/2 < H < 1$, equation (13) is the same as expression (6.4) of Brockwell and Marquardt [12] up to a factor that involves the Hurst parameter only. The asymptotic expression (14) implies that for $0 < H < 1/2$, $\gamma_Y(h) < 0$ when $h \to \infty$, which shows that the CARFIMA$(p, H, q)$ model with $0 < H < 1/2$ is antipersistent. In contrast, the model with $1/2 < H < 1$ is of long-memory type.

The autocovariance function established in Theorem 2 can be used to compute the spectral density function of $\{Y_t, t \in R\}$ stated in Theorem 3 below.

**Theorem 3.** *For $0 < H < 1$, the spectral density function of $\{Y_t, t \in R\}$ is given by*

$$f_Y(w) = \frac{\sigma^2}{2\pi}\Gamma(2H+1)\sin(\pi H)|w|^{1-2H}\frac{|\beta(\mathrm{i}w)|^2}{|\alpha(\mathrm{i}w)|^2}, \qquad \omega \in (-\infty, \infty). \tag{15}$$

For $0 < H < 1/2$, equation (15) implies that $\int_{-\infty}^{\infty} \gamma_Y(\tau)\,\mathrm{d}\tau = 2\pi f_Y(0) = 0$. The equation also shows that the spectral density function of the CARFIMA$(p, H, q)$ process is essentially a product of the spectral density of the ARMA process and the spectral density of the fractional Gaussian noise. Thus, the CARFIMA$(p, H, q)$ model is generated by applying an ARMA filter to the fractional Gaussian noise. Furthermore, the CARFIMA model is antipersistent if the fractional Gaussian noise is antipersistent, whereas it is long-memory if the noise is long-memory. Compared to fractional Gaussian noise, the CARFIMA model displays a much wider spectrum of autocovariance patterns, including non-monotone autocovariance functions.

One major problem with continuous-time modeling is the identifiability of the continuous-time model, given discrete-time data. Let $\{Y_{ih}\}_{i=1,\ldots,N}$ be the observations sampled from a stationary CARFIMA$(p, H, q)$ process, where $h$ is the step size. By the aliasing formula (Priestley [30]), the spectral density of $\{Y_{ih}\}_{i=1,\ldots,N}$ equals

$$f_h(\omega;\theta,\sigma^2) = \frac{1}{h}\sum_{k \in Z} f_Y\left(\frac{\omega + 2k\pi}{h}\right), \qquad \omega \in [-\pi, \pi], \tag{16}$$

where $f_Y(\cdot)$ is as defined in equation (15). Using the frequency domain method, Tsai and Chan [38] showed that the CARFIMA$(p, H, q)$ model with $1/2 < H < 1$ is identifiable, given regularly spaced discrete-time data. Specifically, they showed that, for $(\theta_1, \sigma_1^2) \neq (\theta_2, \sigma_2^2)$, the set $\{\omega | f_h(\omega;\theta_1,\sigma_1^2) \neq f_h(\omega;\theta_2,\sigma_2^2)\}$ has positive Lebesgue measure if $1/2 < H < 1$. Because the spectral density function of $\{Y_{ih}\}_{i=1,\ldots,N}$ is given by the same form as (16) for $0 < H < 1/2$ and $1/2 < H < 1$, it can be shown by similar arguments that the CARFIMA model with $0 < H < 1/2$ is also identifiable (see Section 3). The identifiability problem with $H = 1/2$ is more difficult; see [38] for further discussion. In summary, we have the following theorem on identifiability.

**Theorem 4.** *Let $Y = \{Y_{t_i}\}_{i=1}^N$ be sampled from a stationary (Gaussian) CARFIMA$(p, H, q)$ process given by equation (1), where $0 < H < 1$, $H \neq 1/2$, $\alpha(\cdot)$ and $\beta(\cdot)$ have no common zeros, all roots of $\alpha(z) = 0$ and the roots of $\beta(z) = 0$ have negative real parts. If the step size $t_i = ih$ with $h > 0$, then for $(\theta_1, \sigma_1^2) \neq (\theta_2, \sigma_2^2)$, the set $\{\omega | f_h(\omega;\theta_1,\sigma_1^2) \neq f_h(\omega;\theta_2,\sigma_2^2)\}$ has positive Lebesgue measure.*

We note that the roots of the polynomial $\alpha(\cdot)$ are the same as the eigenvalues of the matrix $A$ and the condition on the roots of $\alpha(z) = 0$ is necessary for the stationarity of the process, whereas the condition on $\beta(z) = 0$ is akin to the invertibility condition for discrete-time processes.

## 3. Proofs

**Proof of equation (9).** Consider simple functions of the form $f(u) = \sum_{i=0}^{m-1} c_i 1_{(s_i, s_{i+1}]}(u)$, where $s_m = s$, and $g(v) = \sum_{j=0}^{n-1} d_j 1_{(t_j, t_{j+1}]}(v)$, where $t_n = t$. Then, the left-hand side of (9) becomes

$$\operatorname{Cov}\left(\int_{-\infty}^{s} f(u)\, dB_u^H, \int_{-\infty}^{t} g(v)\, dB_v^H\right)$$

$$= \operatorname{Cov}\left(\sum_{i=0}^{m-1} c_i(B_{s_{i+1}}^H - B_{s_i}^H), \sum_{j=0}^{n-1} d_j(B_{t_{j+1}}^H - B_{t_j}^H)\right) \tag{17}$$

$$= \tfrac{1}{2}\sum_{i=0}^{m-1}\sum_{j=0}^{n-1} c_i d_j \{|s_{i+1} - t_j|^{2H} + |s_i - t_{j+1}|^{2H} - |t_{j+1} - s_{i+1}|^{2H} - |s_i - t_j|^{2H}\}.$$

The first term of the right-hand side of (9) is

$$H f(s) \int_{-\infty}^{t} |s - v|^{2H-1} \operatorname{sgn}(s - v) g(v)\, dv$$

$$= H c_{m-1} \sum_{j=0}^{n-1} d_j \int_{t_j}^{t_{j+1}} |s_m - v|^{2H-1} \operatorname{sgn}(s_m - v)\, dv \tag{18}$$

$$= \tfrac{1}{2} c_{m-1} \sum_{j=0}^{n-1} d_j \{|s_m - t_j|^{2H} - |s_m - t_{j+1}|^{2H}\}.$$

If we let $c_{-1} = 0$, then the second term of the right-hand side of (9) is

$$H \int_{-\infty}^{t} \int_{-\infty}^{s} g(v) |u - v|^{2H-1} \operatorname{sgn}(v - u)\, df(u)\, dv$$

$$= H \int_{-\infty}^{t} g(v) \sum_{i=0}^{m-1} |s_i - v|^{2H-1} \operatorname{sgn}(v - s_i)(c_i - c_{i-1})\, dv$$

$$= H \sum_{i=0}^{m-1} (c_i - c_{i-1}) \sum_{j=0}^{n-1} d_j \int_{t_j}^{t_{j+1}} |s_i - v|^{2H-1} \operatorname{sgn}(v - s_i)\, dv \tag{19}$$

$$= \tfrac{1}{2} \sum_{i=0}^{m-1}(c_i - c_{i-1}) \sum_{j=0}^{n-1} d_j\{|s_i - t_{j+1}|^{2H} - |s_i - t_j|^{2H}\}$$

$$= \tfrac{1}{2} \sum_{i=0}^{m-1}\sum_{j=0}^{n-1} c_i d_j\{|s_i - t_{j+1}|^{2H} - |s_i - t_j|^{2H}\}$$

$$- \tfrac{1}{2} \sum_{i=0}^{m-2}\sum_{j=0}^{n-1} c_i d_j\{|s_{i+1} - t_{j+1}|^{2H} - |s_{i+1} - t_j|^{2H}\}.$$

Therefore, by (18) and (19), the right-hand side of (9) becomes

$$\tfrac{1}{2} \sum_{i=0}^{m-1}\sum_{j=0}^{n-1} c_i d_j\{|s_i - t_{j+1}|^{2H} - |s_i - t_j|^{2H} + |s_{i+1} - t_j|^{2H} - |s_{i+1} - t_{j+1}|^{2H}\},$$

which is the same as equation (17). This proves the validity of equation (9) for simple functions. □

**Proof of Theorem 1.** The proof for $H = 1/2$ is trivial. For the proof where $1/2 < H < 1$, see Tsai and Chan [37]. We now consider the case where $0 < H < 1/2$. The proof of the first part of the theorem is similar to that of Theorem 1(a) in [37] and is hence omitted. For the proof of equation (10), first note that equations (4) and (6) are essentially equivalent if $\alpha_0 = 0$. Since, by (6), we can write $X_0 = \sigma \int_{-\infty}^0 e^{-Au} \delta_p \, dB_u^H$ and $B_t = \int_0^t dB_u^H$, by equation (9), we have

$$\begin{aligned}
\text{cov}(Y_0, B_t^H) &= \text{cov}\left(\sigma \int_{-\infty}^0 \beta' e^{-Au} \delta_p \, dB_u^H, \int_0^t dB_u^H\right) \\
&= H\sigma \beta' \delta_p \int_0^t |-v|^{2H-1} \text{sgn}(-v) \, dv \qquad (20) \\
&\quad - H\sigma \int_0^t \int_{-\infty}^0 |u-v|^{2H-1} \text{sgn}(v-u) \beta' A e^{-Au} \delta_p \, du \, dv \\
&= -\frac{\sigma}{2} \beta' \delta_p t^{2H} - \frac{\sigma}{2} \int_0^\infty \{(t+u)^{2H} - u^{2H}\} \beta' A e^{Au} \delta_p \, du.
\end{aligned}$$

Now, based on the integration by parts technique, the equality in (10) follows from equation (20). The stationary mean of $\{Y_t\}$ follows from equation (5) and the subsequent discussion. □

**Proof of Theorem 2.** The proofs of 2(a) and (b) with $H = 1/2$ are trivial. First, we prove part (a) with $0 < H < 1/2$. By equations (7), (9) and routine calculus, we have

$$\gamma_Y(h) = \sigma^2 H \int_{-\infty}^{t} \beta' \delta_p \delta_p' e^{A'(t-v)} \beta (t+h-v)^{2H-1} \, dv$$

$$- H\sigma^2 \int_{-\infty}^{t} \int_{0}^{\infty} \beta' A e^{A(t+h+w-v)} \delta_p \beta' e^{A(t-v)} \delta_p w^{2H-1} \, dw \, dv \quad (21)$$

$$+ H\sigma^2 \int_{-\infty}^{t} \int_{0}^{t+h-v} \beta' A e^{A(t+h-w-v)} \delta_p \beta' e^{A(t-v)} \delta_p w^{2H-1} \, dw \, dv.$$

Now, equation (11) follows from equation (21) and equation (6.20) of Karatzas and Shreve [22], namely, $AV^* + V^*A' = -\sigma^2 \delta_p \delta_p'$. For $1/2 < H < 1$, equation (11) follows from Theorem 1(c) of Tsai and Chan [37] and the integration by parts technique.

Below, we prove part 2(b) with $0 < H < 1/2$. The proof with $1/2 < H < 1$ is similar and is hence omitted. When all the eigenvalues of $A$ have negative real parts and are all distinct, Brockwell and Marquardt [12], equation (2.15), show that, for $h \geq 0$,

$$\beta' e^{Ah} \delta_p = \sum_{i=1}^{p} \frac{\beta(\lambda_i)}{\alpha^{(1)}(\lambda_i)} e^{\lambda_i h}. \quad (22)$$

Differentiating the above equation with respect to $h$ on both sides, we have

$$\beta' A e^{Ah} \delta_p = \sum_{i=1}^{p} \frac{\beta(\lambda_i)\lambda_i}{\alpha^{(1)}(\lambda_i)} e^{\lambda_i h}. \quad (23)$$

By expressions (21), (22), (23) and routine calculus, we obtain

$$\gamma_Y(h) \quad (24)$$

$$= -\frac{\sigma^2}{2}\Gamma(2H+1)\left[\sum_{i=1}^{p} \frac{\beta(\lambda_i)}{\alpha^{(1)}(\lambda_i)}\{e^{\lambda_i h}(-\lambda_i)^{1-2H} + e^{\lambda_i h}\lambda_i^{1-2H}P(2H,\lambda_i h)\}\right.$$

$$\times \sum_{j=1}^{p} \frac{\beta(\lambda_j)}{\alpha^{(1)}(\lambda_j)(\lambda_i+\lambda_j)}$$

$$+ \sum_{j=1}^{p} \frac{\beta(\lambda_j)}{\alpha^{(1)}(\lambda_j)} \quad (25)$$

$$\times \{e^{-\lambda_j h}(-\lambda_j)^{1-2H} - e^{-\lambda_j h}(-\lambda_j)^{1-2H}P(2H,-\lambda_j h)\}$$

$$\left.\times \sum_{i=1}^{p} \frac{\beta(\lambda_i)}{\alpha^{(1)}(\lambda_i)(\lambda_i+\lambda_j)}\right].$$

By the partial fraction (Feller [18], page 276)

$$\frac{\beta(s)}{\alpha(s)} = \sum_{j=1}^{p} \frac{\beta(\lambda_j)}{\alpha^{(1)}(\lambda_j)(s-\lambda_j)}, \qquad (26)$$

we have

$$-\frac{\beta(-\lambda_i)}{\alpha(-\lambda_i)} = \sum_{j=1}^{p} \frac{\beta(\lambda_j)}{\alpha^{(1)}(\lambda_j)(\lambda_i+\lambda_j)}. \qquad (27)$$

Equation (12) now follows from equations (25) and (27). This proves (b).

To prove (c) with $0 < H < 1/2$, we note that $P(a,z) = 1 - \Gamma(a,z)/\Gamma(a)$. Furthermore, by equation (6.5.32) in [1], $\Gamma(a,z) \sim z^{a-1}\mathrm{e}^{-z}\{1+(a-1)/z\}$ when $z \to \infty$. Thus, as $h \to \infty$, we have

$$u(H,\lambda,h) \sim -\frac{4H(2H-1)}{\Gamma(2H+1)\lambda}h^{2H-2}. \qquad (28)$$

Hence, by equations (12) and (28), as $h \to \infty$, we have

$$\gamma_Y(h) \sim -2\sigma^2 H(2H-1)h^{2H-2} \sum_{i=1}^{p} \frac{\beta(\lambda_i)\beta(-\lambda_i)}{\alpha^{(1)}(\lambda_i)\alpha(-\lambda_i)\lambda_i}. \qquad (29)$$

Now, $1/(\lambda_i\lambda_j) = 1/\{\lambda_i(\lambda_i+\lambda_j)\} + 1/\{\lambda_j(\lambda_i+\lambda_j)\}$ and expression (27) implies that

$$\begin{aligned}
\sum_{i=1}^{p} & \frac{\beta(\lambda_i)\beta(-\lambda_i)}{\alpha^{(1)}(\lambda_i)\alpha(-\lambda_i)\lambda_i} \\
&= -\sum_{i=1}^{p}\sum_{j=1}^{p} \frac{\beta(\lambda_i)\beta(\lambda_j)}{\alpha^{(1)}(\lambda_i)\alpha^{(1)}(\lambda_j)\lambda_i(\lambda_i+\lambda_j)} \\
&= -\frac{1}{2}\sum_{i=1}^{p}\sum_{j=1}^{p} \frac{\beta(\lambda_i)\beta(\lambda_j)}{\alpha^{(1)}(\lambda_i)\alpha^{(1)}(\lambda_j)\lambda_i(\lambda_i+\lambda_j)} \\
&\quad -\frac{1}{2}\sum_{i=1}^{p}\sum_{j=1}^{p} \frac{\beta(\lambda_i)\beta(\lambda_j)}{\alpha^{(1)}(\lambda_i)\alpha^{(1)}(\lambda_j)\lambda_j(\lambda_i+\lambda_j)} \\
&= -\frac{1}{2}\sum_{i=1}^{p}\sum_{j=1}^{p} \frac{\beta(\lambda_i)\beta(\lambda_j)}{\alpha^{(1)}(\lambda_i)\alpha^{(1)}(\lambda_j)\lambda_i\lambda_j} \\
&= -\frac{1}{2}\frac{\beta^2(0)}{\alpha^2(0)},
\end{aligned} \qquad (30)$$

where the last equality follows from the identity (26). The proof of (c) with $0 < H < 1/2$ now follows from expressions (29) and (31). Note that the proof of (c) with $1/2 < H < 1$ is similar to that of $0 < H < 1/2$ and is therefore omitted. □

**Proof of Theorem 3.** For the proof with $1/2 < H < 1$, readers can refer to Tsai and Chan [38]. For case where $H = 1/2$ is trivial, the proof with $0 < H < 1/2$ is similar to that of Theorem 1 in [38] and is hence omitted. □

**Proof of Theorem 4.** For the proof with $1/2 < H < 1$, see Tsai and Chan [38]. The following proof with $0 < H < 1/2$ is similar to that in [38].

If we write $f_c(\omega) = f_Y(\omega) = g(H)L(\omega)|\omega|^{1-2H}$, where $g(H) = \Gamma(2H+1)\sin(\pi H)/(2\pi)$, $L(\omega) = \sigma^2|\beta(i\omega)|^2/|\alpha(i\omega)|^2$, then $f_c^{(1)}(\omega) = g(H)|\omega|^{-2H}\{|\omega|L^{(1)}(\omega) + (1-2H)L(\omega)\}$ and

$$\log(f^{(1)}(\omega)) = \log(f_c^{(1)}(\omega)) + \log\left\{1 + \frac{\sum_{k \neq 0} f_c^{(1)}(\omega + 2k\pi)}{f_c^{(1)}(\omega)}\right\}$$

$$= \log g(H) - 2H \log|\omega| + \log\{|\omega|L^{(1)}(\omega) + (1-2H)L(\omega)\} \qquad (31)$$

$$+ \log(1 + R(\omega)),$$

where $R(\omega) = R_2(\omega)/R_1(\omega)$, $R_1(\omega) = f_c^{(1)}(\omega)$ and $R_2(\omega) = \sum_{k \neq 0} f_c^{(1)}(\omega + 2k\pi)$. It follows from $\lim_{\omega \to 0}\{\log(f^{(1)}(\omega))/\log|\omega|\} = -2H$ that the Hurst parameter $H$ is identifiable. Following Tsai and Chan [38], it suffices to show that, given $f(\omega), \omega \in [-\pi, \pi]$, we can determine $L$ and all of its higher derivatives at $\omega = 0$. Below, we show that $L(0)$ and $L^{(1)}(0)$ are identifiable. The identifiability of $L^{(k)}(0)$ for $k \geq 2$ can be proven similarly by the arguments in [38]. The identifiability of $L(0)$ simply follows from the fact that

$$\lim_{\omega \to 0}\{\log(f^{(1)}(\omega)) - \log g(H) + 2H\log|\omega| - \log(1-2H)\}$$

$$= \lim_{\omega \to 0}\left[\log L(\omega) + \log\left\{1 + \frac{|\omega|L^{(1)}(\omega)}{(1-2H)L(\omega)}\right\} + \log(1+R(\omega))\right]$$

$$= \log L(0).$$

Next, we prove the identifiability of $R_2(0)$, which is needed to prove the identifiability of $L^{(1)}(0)$. By equation (31), we have

$$\frac{\partial}{\partial \omega}\{\log(f^{(1)}(\omega)) + 2H\log|\omega|\}$$

$$= \frac{|\omega|L^{(2)}(\omega) + (2-2H)L^{(1)}(\omega)}{|\omega|L^{(1)}(\omega) + (1-2H)L(\omega)} + \frac{R^{(1)}(\omega)}{1+R(\omega)}, \qquad (32)$$

where

$$\frac{R^{(1)}(\omega)}{1+R(\omega)} = \frac{R_2^{(1)}(\omega)}{(1+R(\omega))R_1(\omega)} - \left\{\frac{R_1^2(\omega)}{R_2(\omega)R_1^{(1)}(\omega)} + \frac{R_1(\omega)}{R_1^{(1)}(\omega)}\right\}^{-1}. \qquad (33)$$

By equations (32) and (33), we have

$$\lim_{\omega \to 0} \left[ |\omega|^{1-2H} \frac{\partial}{\partial \omega} \{\log(f^{(1)}(\omega)) + 2H \log |\omega|\} \right]$$

$$= \lim_{\omega \to 0} \left\{ \frac{|\omega|^{1-2H} R^{(1)}(\omega)}{1 + R(\omega)} \right\}$$

$$= \lim_{\omega \to 0} \left\{ -\frac{|\omega|^{1-2H} R_2(\omega) R_1^{(1)}(\omega)}{R_1^2(\omega)} \right\}$$

$$= \lim_{\omega \to 0} \left[ -\frac{R_2(\omega)\{|\omega|^2 L^{(2)}(\omega) + 2(1-2H)|\omega|L^{(1)}(\omega) - 2H(1-2H)L(\omega)\}}{g(H)\{|\omega|L^{(1)}(\omega) + (1-2H)L(\omega)\}^2} \right]$$

$$= \frac{2HR_2(0)}{g(H)(1-2H)L(0)},$$

hence $R_2(0)$ is identifiable. Now, equation (33) implies that

$$\lim_{\omega \to 0} \left\{ \frac{R^{(1)}(\omega)}{1 + R(\omega)} - \frac{2HR_2(0)|\omega|^{2H-1}}{g(H)(1-2H)L(0)} \right\}$$

$$= \lim_{\omega \to 0} \left\{ -\frac{R_2(\omega) R_1^{(1)}(\omega)}{R_1^2(\omega)} - \frac{2HR_2(0)|\omega|^{2H-1}}{g(H)(1-2H)L(0)} \right\}$$

$$= \lim_{\omega \to 0} \left[ -\frac{R_2(\omega)|\omega|^{2H-1}\{|\omega|^2 L^{(2)}(\omega) + 2(1-2H)|\omega|L^{(1)}(\omega) - 2H(1-2H)L(\omega)\}}{g(H)\{|\omega|L^{(1)}(\omega) + (1-2H)L(\omega)\}^2} \right. \tag{34}$$

$$\left. - \frac{2HR_2(0)|\omega|^{2H-1}}{g(H)(1-2H)L(0)} \right]$$

$$= \frac{2H}{g(H)} \lim_{\omega \to 0} \left[ |\omega|^{2H-1} \left\{ \frac{(1-2H)L(\omega)R_2(\omega)}{\{|\omega|L^{(1)}(\omega) + (1-2H)L(\omega)\}^2} - \frac{R_2(0)}{(1-2H)L(0)} \right\} \right]$$

$$= 0$$

because $R_2(\omega) = R_2(0) + \omega O(1)$ and $L(\omega) = L(0) + \omega O(1)$ for $\omega$ tending to 0. Therefore, by (32), (33) and (34),

$$\lim_{\omega \to 0} \left[ \frac{\partial}{\partial \omega} \{\log(f^{(1)}(\omega)) + 2H \log |\omega|\} - \frac{2HR_2(0)|\omega|^{2H-1}}{g(H)(1-2H)L(0)} \right]$$

$$= \lim_{\omega \to 0} \frac{|\omega|L^{(2)}(\omega) + (2-2H)L^{(1)}(\omega)}{|\omega|L^{(1)}(\omega) + (1-2H)L(\omega)}$$

$$= \frac{(2-2H)L^{(1)}(0)}{(1-2H)L(0)},$$

which proves the identifiability of $L^{(1)}(0)$. □

## 4. Conclusion

We have proposed a unified continuous-time framework that is useful for studying time series with short memory, long memory and antipersistence. The identifiability of the CARFIMA process with discrete-time data established in Theorem 4 is a fundamental feature that makes the model practical for data analysis. Therefore, in future research, it would be both interesting and useful to study the statistical inference of antipersistent CARFIMA models for data analysis.

## Acknowledgements

We wish to thank Kung-Sik Chan, Paul Dunne and two referees for their careful reading and valuable suggestions leading to a streamlined version of the paper. We also wish to thank Academia Sinica, the National Science Council (NSC 91-2118-M-001-011), R.O.C. and the National Science Foundation (DMS-0405267) for their support.